\documentclass[11pt,reqno, a4paper]{amsart}
\usepackage{amsfonts}
\usepackage{amssymb}
\usepackage{txfonts}
\usepackage{mathrsfs}
\usepackage{bbm}
\usepackage{enumerate}

\newtheorem{thm}{Theorem}[section]
 
 \newtheorem{lem}[thm]{Lemma}
 \newtheorem{prop}[thm]{Proposition}
 
 \newtheorem{rem}[thm]{Remark}

 \newcommand{\vpi}{\varphi}
 \newcommand{\vstar}{v_{*}}

 \newcommand{\p}{\partial}
 
 \newcommand{\s}{\mathcal{S}}
 \newcommand{\Real}{\mathbb{R}}
 
 \newcommand{\Natural}{\mathbb{N}}
 
 \newcommand{\norm}[1]{\Vert#1\Vert}
 
 \newcommand{\abs}[1]{\left\vert#1\right\vert}
 \newcommand{\bigabs}[1]{\big\vert#1\big\vert}
 \newcommand{\set}[1]{\left\{#1\right\}}
 \newcommand{\bigset}[1]{\big\{#1\big\}}
 \newcommand{\inner}[1]{\left(#1\right)}
 \newcommand{\biginner}[1]{\big(#1\big)}
 
 \newcommand{\com}[1]{\big[#1\big]}
 \newcommand{\reff}[1]{(\ref{#1})}
 \newcommand{\V}{\,\,|\,\,\,}

\begin{document}
 \title{Gevrey regularity for the solution of the spatially homogeneous
Landau equation}

\author{Hua Chen, Wei-Xi Li \and Chao-Jiang Xu}


 \address{School of Mathematics and Statistics, Wuhan University, Wuhan 430072, P. R. China}
 \email{chenhua@whu.edu.cn}

\address{School of Mathematics and Statistics, Wuhan University,Wuhan 430072, P. R. China}
 \email{weixi.li@yahoo.com}

 \address{School of Mathematics and Statistics, Wuhan University, Wuhan 430072, P. R. China \newline
 \indent and \newline \indent Universit\'{e} de Rouen, UMR 6085-CNRS, Math\'{e}matiques, Avenue de
l'Universit\'{e}, BR.12, \newline \indent F76801 Saint Etienne du
Rouvray, France}
 \email{Chao-Jiang.Xu@univ-rouen.fr}

\begin{abstract}
  In this paper, we study the
Gevrey  class regularity for solutions of the spatially homogeneous
Landau equations in the hard potential case and the Maxwellian
molecules case.
\end{abstract}
\maketitle

\section{Introduction}

In this paper we study the smoothness effects of solutions for the
following Cauchy problem of the spatially homogeneous Landau
equation
\begin{equation}\label{Landau}
\left\{
\begin{array}{ll}
   \partial_t f=\nabla_{v}\cdot\bigset{\int_{{\Real}^3}
   a(v-\vstar)[f(\vstar)\nabla_v
   f(v)-f(v)\nabla_{v}f(\vstar)]d\vstar},\\
  f(0,v)=f_0(v),
\end{array}
\right.
\end{equation}
where $f(t,v)\geq 0$ stands for the density of particles with
velocity $v\in{\Real}^3$ at time $t\geq0$, and $(a_{ij})$ is a
nonnegative symmetric matrix given by
\begin{equation}\label{coe}
  a_{ij}(v)=\inner{\delta_{ij}-\frac{v_iv_j}{\abs v^2}}\abs
  v^{\gamma+2}.
\end{equation}
We only consider here the condition $\gamma\in[0, 1].$ It's called
the hard potential case when $\gamma\in]0, 1]$ and the Maxwellian
molecules case when $\gamma=0$.

Set $c=\sum_{i,j=1}^3\partial_{v_iv_j}a_{ij}=-2(\gamma+3)\abs
 v^\gamma$ and
 \[
 \bar a_{ij}(t,v)=\inner{a_{ij}*f}(t,v)=\int_{{\Real}^3}a_{ij}(v-\vstar)f(t,\vstar)d\vstar,
 \quad \bar c=c*f.
 \]
Then the Cauchy problem \reff{Landau} can be rewritten as the
following form,
\begin{equation}\label{Landau'}
\left\{
  \begin{array}{ll}
  \partial_t f=\sum_{i,j=1}^3
  \bar a_{ij}\partial_{v_iv_j}f-\bar cf,\\
  f(0,v)=f_0(v).
  \end{array}
\right.
\end{equation}
This is a non-linear  diffusion equation,  and the coefficients
$\bar a_{ij}, \bar c$ depend on the solution $f.$

Here we are mainly concerned with the Gevrey class regularity for
the solution of the Landau equation. This equation is obtained as a
limit of the Boltzmann equation when the collisions become grazing
(see \cite{des1} and references therein). Recently, a lot of
progress has been made on the study of the Sobolev regularizing
property, cf. \cite{C-D-H, desv-vill1, Guo, villani, villani2} and
references therein, which shows that in some sense the Landau
equation can be regarded as a non-linear and non-local analog of the
hypo-elliptic Fokker-Planck equation. That means the weak solution,
which constructed under rather weak hypothesis on the initial datum,
will become smooth or, even more, rapidly decreasing in $v$ at
infinity. This behavior is quite similar to that of the spatially
homogeneous Boltzmann equation without cut-off (see \cite{al-2,
desv-wen1} for more details).

In the Gevrey class frame, some results have been obtained
concerning the propagation property for  solutions of the Boltzmann
equation, e.g. the solutions having the Gevrey regularity have been
constructed in \cite{ukai} for initial data having the same Gevrey
regularity, and the uniform propagation in all time of the Gevrey
regularity has been proved in \cite{des-fur-ter} in the case of
Maxwellian molecules. Recently a general Gevrey regularity result
have been given in \cite{MUXY1} for spatially homogeneous and linear
Boltzmann equation for any initial data. On the other hand, the
local Gevrey regularity for all variables $t,x,v$ is obtained in
\cite {chen-li-xu} for some semi-linear Fokker-Planck equations,
which implies that, in this case, there are also the smoothness
effects which is similar to the heat equation case.

Now we give some notations used throughout the paper.  For a
multi-index $\alpha=(\alpha_1,\alpha_2,\alpha_3),$ we denote
$\abs\alpha=\alpha_1+\alpha_2+\alpha_3$,
$\alpha!=\alpha_1!\alpha_2!\alpha_3!$ and
$\p^\alpha=\p^{\alpha_1}_{v_1}\p^{\alpha_2}_{v_2}\p^{\alpha_3}_{v_3}$.
We say $\beta=(\beta_1,\beta_2, \beta_3)\leq (\alpha_1,\alpha_2,
\alpha_3)=\alpha$ if $\beta_i\leq \alpha_i$ for each $i$. Denote by
$M(f(t))$, $E(f(t))$ and $H(f(t))$ respectively  the mass, energy
and entropy of the function $f(t)$, i.e.,
\begin{eqnarray*}
M(f(t))=\int_{{\Real}^3}f(t,v)\,dv,\quad
E(f(t))={1\over2}\int_{{\Real}^3}f(t,v)\abs v^2\,dv,
\end{eqnarray*}
\[H(f(t))=\int_{{\Real}^3}f(t,v)\log f(t,v)\,dv.\]
and denote $M_0=M(f(0))$, $E_0=E(f(0))$ and $H_0=H(f(0)).$ We know
that the solution of the Landau equation satisfies the formal
conservation laws:
\[M(f(t))=M_0 ,\quad
E(f(t))=E_0,\quad H(f(t))\leq H_0,\qquad \forall\:t\geq 0.\] Also we
adopt the following notations,
\[
 \norm {\p^\alpha f(t,\cdot)}_{L_s^p}^p=\norm {\p^\alpha f(t)}_{L_s^p}^p=\int_{{\Real}^3}\abs{\p^\alpha f(t,v)}^p\inner{1+\abs
 v^2}^{s/2}dv, \mbox{ for }p\geq 1,
\]
\[
 \norm{f(t,\cdot)}_{H^m_s}^2= \norm{f(t)}_{H^m_s}^2
 =\sum_{\abs\alpha\leq m}\norm{\p^\alpha f(t,\cdot)}_{L^2_s}^2.\]

Before stating our main result, let us recall the definition of the
Gevrey class function space $G^\sigma({\Real}^N),$ where
$\sigma\geq1$ is the Gevrey index (cf. \cite{CR, Ro}). Let $u$ be a
real function defined in ${\Real}^N.$ We say $u\in
G^\sigma({\Real}^N)$ if $u\in C^\infty({\Real}^N)$ and  there exists
a positive constant $C$ such that for all multi-indices
$\alpha\in{\Natural}^N,$ we have
\[
  \norm{\partial^\alpha{u}}_{L^2({\Real}^3)}\leq C^{|\alpha|+1}
  (\abs\alpha!)^\sigma.
\]
We denote by $G_0^\sigma({\Real}^3)$ the space of Gevrey functions
with compact support. Note that $G^1({\Real}^N)$ is the space of
real analytic functions.

In the hard potential case, the existence, uniqueness and Sobolev
regularity of the weak solution  had been studied by Desvillettes
and Villani (see \cite{desv-vill1}, Theorem 5, Theorem 6 and Theorem
7). Actually they proved that, under rather weak assumptions on the
initial datum (for instance if $f_0\in L^1_{2+\delta}$ with
$\delta>0$), there exists a weak solution $f$ of the Cauchy problem
\reff{Landau'} such that for all time $t_0>0,$ all integer $m\geq
0,$ and all $s>0,$
\[\sup\limits_{t\geq t_0}\norm{f(t,\cdot)}_{H^m_s}\leq C,\]
where $C$ is a constant depending only on $\gamma$, $M_0$, $E_0$,
$H_0$, $m,s$ and $t_0.$ Furthermore they proved that $f(t,v)\in
C^{\infty}\biginner{{\Real}_t^+; \mathcal {S}({\Real}_v^3)},$ where
${\Real}_t^+=]0, +\infty[$ and $\s ({\Real}_v^3)$ denotes the space
of smooth functions rapidly decreasing at infinity.   If $f_0\in
L^2_p$ with $p>5\gamma+15,$ then the Cauchy problem admits a unique
smooth solution .

In the Maxwellian case, Villani \cite{villani2} proved that the
Cauchy problem admits a unique classical solution $f$ for any
initial datum and for all $t>0,$ and that $f$ is bounded and lies in
$C^\infty({\Real}_v^3).$

\smallskip
Now starting from  the smooth solution,  we state our main result on
the Gevrey regularity as follows.

\begin{thm}\label{main result}
Let $f_0$ be the initial datum with finite mass, energy and entropy
and $f$ be  any solution of the Cauchy problem \reff{Landau'} such
that for all $t_0$, $t_1$ with $0<t_0<t_1<+\infty,$ and all integer
$m\geq0,$
\begin{equation}\label{smooth}
 \sup\limits_{t\in[t_0, t_1]}\norm {f(t,\cdot)}_{H^m_\gamma}<+\infty.
\end{equation} Then for any number $\sigma>1,$  we have
$f(t, \cdot)\in G^\sigma({\Real}^3)$ for all time $t>0$.
\end{thm}

\begin{rem}
Our result, which is given here for the space dimension to be equal
to $3$, will also be true for any space dimension.
\end{rem}

This paper includes three sections. The proof of the main result
Theorem \ref{main result} will be given in Section \ref{section2},
and in Section \ref{section3} we shall mainly prove Lemma
\ref{crucial}, which is crucial in the proof of Section
\ref{section2}.

\section{Proof of the main results}
\label{section2}

This section is devoted to the proof of  the main result. In the
sequel, we always use $\sum_{1\leq \abs\beta\leq \abs\mu}$ to
denotes the summation over all the multi-indices $\beta$ satisfying
$\beta\leq\mu$ and $1\leq \abs\beta\leq\abs\mu$. Likewise
$\sum_{1\leq \abs\beta\leq \abs\mu-1}$ denotes the summation over
all the multi-indices $\beta$ satisfying $\beta\leq\mu$ and $1\leq
\abs\beta\leq\abs\mu-1.$ Firstly we have

\begin{lem}
For any $\sigma>1,$ there exists a constant $C_\sigma,$ depending
only on $\sigma,$ such that for all multi-indices
$\mu\in{\Natural}^3, \abs\mu\geq 1,$
\begin{equation}\label{sum1}
  \sum_{1\leq \abs\beta\leq \abs\mu}\frac{1}{\abs\beta^{3}}\leq
  C_\sigma\abs\mu^{\sigma-1},
\end{equation}
and
\begin{equation}\label{sum2}
  \sum_{1\leq \abs\beta\leq\abs\mu-1}\frac{1}{\abs\beta^2(\abs\mu-\abs\beta)}\leq
  C_\sigma\abs\mu^{\sigma-1}.
\end{equation}
\end{lem}

{\bf Proof~} For each positive integer $l,$ we denote by
$N\set{\abs\beta=l}$ the number of the multi-indices $\beta$ with
$\abs\beta=l.$ In the case when the space dimension equals to 3, one
has
\[N\set{\abs\beta=l}=\frac{(l+2)!}{2!\,l!}={1\over2}(l+1)(l+2).\]
It's easy to see that
\begin{eqnarray*}
  \sum_{1\leq \abs\beta\leq
  \abs\mu}\frac{1}{\abs\beta^{3}}\leq\sum_{l=1}^{
  \abs\mu}\sum_{\abs\beta=l}\frac{1}{l^{3}}=\sum_{l=1}^{
  \abs\mu}\frac{N\set{\abs\beta=l}}{l^{3}}
\end{eqnarray*}
We combine these estimates to compute
\begin{eqnarray*}
  \sum_{1\leq \abs\beta\leq
  \abs\mu}\frac{1}{\abs\beta^{3}}\leq{1\over2}\sum_{l=1}^{
  \abs\mu}\frac{(l+1)(l+2)}{l^{3}}\leq 3\sum_{l=1}^{
  \abs\mu}\frac{1}{l}.
\end{eqnarray*}
Together  with the fact that $3\sum_{l=1}^{
  \abs\mu}l^{-1}\leq C_\sigma\abs\mu^{\sigma-1}$  for some constant $C_\sigma,$
 this gives the estimate
\reff{sum1}. In a similar way as above we can prove the estimate
\reff{sum2}. The proof is completed.

The following lemma, which will be proved in the next section, is of
great of use to us. For simplicity, in the following discussions, we
shall use the notation $\gamma-j,$ with $\gamma$ a multi-index and
$j$ an integer, to denotes some multi-index $\tilde\gamma$
satisfying $\tilde\gamma\leq\gamma$ and
$\abs{\tilde\gamma}=\abs\mu-j.$

\begin{lem}\label{crucial}
Let $\sigma>1$. There exist constants $B, C_1, C_2>0$  with $B$
depending only on the Gevrey index $\sigma$  and $C_1, C_2$
depending only on $M_0, E_0, H_0, \sigma$ and $\gamma,$  such that
for all multi-indices $\mu\in{\Natural}^3$ with $\abs\mu\geq2$ and
all $t>0,$ we have
\begin{eqnarray*}
 \frac{d}{dt}&\norm&{\p^\mu f(t)}_{L^2}^2+C_1\norm{\nabla_v\p^\mu  f(t) }_{L^2_\gamma}^2
  \leq C_2\abs\mu^2\norm{\nabla_v\p^{\mu-1} f(t)}_{L^2_\gamma}^2\\
  &+&C_2 \sum_{2\leq\abs\beta\leq\abs\mu}
  C_{\mu }^\beta \norm{\nabla_v\p^{\mu -\beta+1} f(t)}_{L^2_\gamma}
  \cdot\norm{\nabla_v\p^{\mu-1}
  f(t)}_{L^2_\gamma}\cdot\com{G_\sigma(f(t))}_{\beta-2}\\
  &+&C_2\sum_{0\leq\abs\beta\leq\abs\mu}
   C_{\mu }^\beta\norm{
  \p^\beta f(t)}_{L^2_\gamma}\cdot\norm{\nabla_v\p^{\mu-1} f(t)}_{L^2_\gamma}
  \cdot\com{G_\sigma(f(t))}_{{\mu}-\beta},
\end{eqnarray*}
where $C_\mu^\beta=\frac{\mu!}{(\mu-\beta)!\beta!}$ is the binomial
coefficients, and
$$\com{G_\sigma(f(t))}_{\beta-2}=\norm{\p^{\beta-2}f(t)}_{L^2}+B^{\abs\beta-2}\inner{(\abs\beta-2)!}^\sigma.$$

\end{lem}

From now on, $\Omega$ will be used to denote an arbitrary fixed
interval $[T_0, T_1]$ with $0<T_0<T_1<T_0+1.$  We denote
\[\Omega_\rho=[T_0+\rho,T_1-\rho],\quad 0<\rho<{{T_1-T_0}\over2}<{1\over2}.\]
For any integer $k$ with $k\geq 2$ and any $\rho$ with
$0<\rho<(T_1-T_0)/2,$ take a function $\varphi_{\rho,k}(t)\in
C_0^\infty({\Real})$ satisfying $0\leq\vpi_{\rho,k}\leq1$, and Supp
$ \varphi_{\rho,k}\subset\Omega_{\tilde\rho}$ with $
\tilde\rho=\frac{k-1}{k}\rho,$ and  $\varphi_{\rho,k}=1$ in
$\Omega_{\rho}$. it is easy to verify
  \begin{equation}\label{vpi}  \sup\abs{\frac{d^j\varphi_{\rho,k}}{dt^j}}\leq
      \tilde C_j\: (k/\rho)^{j},\indent \forall~j\in\Natural.
  \end{equation}
And for $\tilde\rho=(k-1)\rho/k, ~k\geq2, $  the following simple
fact is clear,
\begin{equation}\label{fact}
 {1\over{\rho}^{k}}\leq{1\over{\tilde\rho}^{k}}={1\over{\rho}^{k}}\times\big({k\over{k-1}}\big)^{k}
\leq  {5\over{\rho}^{k}}.
\end{equation}

Now we are prepare to prove the main results, which can be deduced
easily from the following

\begin{prop}\label{prp}

Let $f_0$ be the initial datum with finite mass, energy and entropy
and $f$ be  any solution of the Cauchy problem satisfying
\reff{smooth}. Then for $\sigma>1$, there exists a constant $A$,
depending only on $T_0, T_1, M_0, E_0, H_0, \gamma$ and $\sigma$,
such that for any $k\in\Natural,$ $k\geq 0,$
\begin{eqnarray*}
 (Q)_k \quad \sup_{t\in\Omega_\rho}\norm{\p^\alpha f(t)}_{L^2}
 +\set{\int_{T_0+\rho}^{T_1-\rho}\|\nabla_v \p^{\tilde\alpha}f(t)\|_{L^2_{\gamma}}^2\,dt}^{1/2}\leq
 \frac{A^{k}}{\rho^k}\com{(k-1)!}^\sigma
\end{eqnarray*}
holds for any  multi-indices $\alpha,\tilde\alpha$ with
$\abs\alpha=\abs{\tilde\alpha}=k$ and all $\rho$ with
$0<\rho<(T_1-T_0)/2.$ Here we assume $(-1)!=0!=1.$
\end{prop}

{\bf Proof~} We use induction on $k$. $(Q)_0, (Q)_1$ obviously hold
if we take $A$ large enough such that
\begin{equation}\label{A}
 A\geq\sup_{s\in[T_0,T_1]}\norm {f(s)}_{H^2_\gamma}+1
\end{equation}
The term on the right hand of \reff{A} is finite by virtue of
\reff{smooth}. Now assuming $(Q)_{k-1}$ holds, we shall show the
truth of $(Q)_k, k\geq 2.$ In this proof $C_j, j\geq 3$, are used to
denote different constants depending only on $T_0, T_1, M_0, E_0,
H_0, \gamma$ and $\sigma$.

Firstly we shall prove
\begin{equation}\label{first term}
\sup_{t\in\Omega_\rho}\norm{\p^\alpha f(t)}_{L^2}\leq {1\over
2}\frac{A^{\abs\alpha}}{\rho^{\abs\alpha}}\com{(\abs\alpha-1)!}^\sigma,\quad
\forall~\abs\alpha=k,\quad \forall~0<\rho<\frac{T_1-T_2}{2}.
\end{equation}
In the following discussion, let $\alpha$ be any multi-index with
$\abs\alpha=k$ and $\rho$ be any number with $0<\rho<(T_1-T_2)/2.$
Applying Lemma \ref{crucial} with $\mu=\alpha,$ we obtain
\begin{eqnarray*}
  &\frac{d}{dt}&\norm{\p^\alpha f(t)}_{L^2}^2+C_1\norm{\nabla_v\p^\alpha  f(t) }_{L^2_\gamma}^2
  \leq C_2\abs\alpha^2\norm{\nabla_v\p^{\alpha-1} f(t)}_{L^2_\gamma}^2\\
  &\indent+&C_2 \sum_{2\leq\abs\beta\leq\abs\alpha}
  C_{\alpha }^\beta \norm{\nabla_v\p^{\alpha -\beta+1} f(t)}_{L^2_\gamma}
  \cdot\norm{\nabla_v\p^{\alpha-1}
  f(t)}_{L^2_\gamma}\cdot\com{G_\sigma(f(t))}_{\beta-2}\\
  &\indent+&C_2\sum_{0\leq\abs\beta\leq\abs\alpha}
   C_{\alpha }^\beta\norm{
  \p^\beta f(t)}_{L^2_\gamma}\cdot\norm{\nabla_v\p^{\alpha-1} f(t)}_{L^2_\gamma}
  \cdot\com{G_\sigma(f(t))}_{{\alpha}-\beta}.
\end{eqnarray*}
Write the last term of the right side  as
\begin{eqnarray*}
  &C_2\norm{
  f(t)}_{L^2_\gamma}\cdot\norm{\nabla_v\p^{\alpha-1} f(t)}_{L^2_\gamma}
  \cdot\com{G_\sigma(f(t))}_{{\alpha}}\\
  &+C_2\sum_{1\leq\abs\beta\leq\abs\alpha}
   C_{\alpha }^\beta\norm{
  \p^\beta f(t)}_{L^2_\gamma}\cdot\norm{\nabla_v\p^{\alpha-1} f(t)}_{L^2_\gamma}
  \cdot\com{G_\sigma(f(t))}_{{\alpha}-\beta}.
\end{eqnarray*}
And then we get
\begin{eqnarray*}
  &\frac{d}{dt}&\norm{\p^\alpha f(t)}_{L^2}^2+C_1\norm{\nabla_v\p^\alpha  f(t) }_{L^2_\gamma}^2
  \leq C_2\abs\alpha^2\norm{\nabla_v\p^{\alpha-1} f(t)}_{L^2_\gamma}^2\\
  &\indent+&C_2\norm{
  f(t)}_{L^2_\gamma}\cdot\norm{\nabla_v\p^{\alpha-1} f(t)}_{L^2_\gamma}
  \cdot\com{G_\sigma(f(t))}_{{\alpha}}\\
  &\indent+&C_2 \sum_{2\leq\abs\beta\leq\abs\alpha}
  C_{\alpha }^\beta \norm{\nabla_v\p^{\alpha -\beta+1} f(t)}_{L^2_\gamma}
  \cdot\norm{\nabla_v\p^{\alpha-1}
  f(t)}_{L^2_\gamma}\cdot\com{G_\sigma(f(t))}_{\beta-2}\\
  &\indent+&C_2\sum_{1\leq\abs\beta\leq\abs\alpha}
   C_{\alpha }^\beta\norm{
  \p^\beta f(t)}_{L^2_\gamma}\cdot\norm{\nabla_v\p^{\alpha-1} f(t)}_{L^2_\gamma}
  \cdot\com{G_\sigma(f(t))}_{{\alpha}-\beta}.
\end{eqnarray*}
Multiplying by $\vpi_{\rho, k}(t)$ the both sides of the above
inequality, one has
\begin{eqnarray*}
  &\frac{d}{dt}&\com{\vpi_{\rho,k}(t)\norm{\p^\alpha f(t)}_{L^2}^2}+C_1\vpi_{\rho,k}(t)\norm{\nabla_v\p^\alpha  f(t)
  }_{L^2_\gamma}^2\\
  &\leq& \frac{d\vpi_{\rho,k}}{dt}\cdot\norm{\p^\alpha f(t)}_{L^2}^2+C_2\cdot\vpi_{\rho,k}(t)\abs\alpha^2\norm{\nabla_v\p^{\alpha-1} f(t)}_{L^2_\gamma}^2\\
  &\quad+&C_2\cdot\vpi_{\rho,k}(t)~\norm{
  f(t)}_{L^2_\gamma}\cdot\norm{\nabla_v\p^{\alpha-1} f(t)}_{L^2_\gamma}
  \cdot\com{G_\sigma(f(t))}_{{\alpha}}\\
  &\quad+&C_2 \cdot\vpi_{\rho,k}(t)\sum_{2\leq\abs\beta\leq\abs\alpha}
  C_{\alpha }^\beta \norm{\nabla_v\p^{\alpha -\beta+1} f(t)}_{L^2_\gamma}
  \cdot\norm{\nabla_v\p^{\alpha-1}
  f(t)}_{L^2_\gamma}\cdot\com{G_\sigma(f(t))}_{\beta-2}\\
  &\quad+&C_2\cdot\vpi_{\rho,k}(t)\sum_{1\leq\abs\beta\leq\abs\alpha}
   C_{\alpha }^\beta\norm{
  \p^\beta f(t)}_{L^2_\gamma}\cdot\norm{\nabla_v\p^{\alpha-1} f(t)}_{L^2_\gamma}
  \cdot\com{G_\sigma(f(t))}_{{\alpha}-\beta}.
\end{eqnarray*}
To simplify the notation, we set
\[\com{G_\sigma(f)}_{\rho,\beta}=\sup_{s\in\Omega_\rho}\com{G_\sigma(f(s))}_{\beta}
=\sup_{s\in\Omega_\rho}\norm{\p^\beta
f(s)}_{L^2}+B^{\abs\beta}(\abs\beta!)^\sigma.\] Recall Supp
$\vpi_{\rho,k}\subset\Omega_{\tilde\rho}$ with
$\tilde\rho=(k-1)\rho/k$ and $\vpi_{\rho,k}(t)=1$ for all
$t\in\Omega_\rho$ and $\vpi_{\rho,k}(T_0)=0.$ Then  for any
$t\in\Omega_\rho,$  we integrate the above inequality over the
interval $[T_0, t]\subset[T_0,T_1-\rho]$ and then use  Cauchy
inequality to get
\begin{eqnarray*}
  &\norm&{\p^\alpha f(t)}_{L^2}^2=\vpi_{\rho,k}(t)\norm{\p^\alpha f(t)}_{L^2}^2-\vpi_{\rho,k}(T_0)~\norm{\p^\alpha
  f(T_0)}_{L^2}^2\\
  &\leq& \sup\abs{\frac{d\vpi_{\rho,k}}{dt}}\int_{T_0+\tilde\rho}^{T_1-\tilde\rho}\norm{\p^{\alpha} f(s)}_{L^2}^2ds+
  C_2\abs\alpha^2\int_{T_0+\tilde\rho}^{T_1-\tilde\rho}\norm{\nabla_v\p^{\alpha-1} f(s)}_{L^2_\gamma}^2ds\\
  &+&C_2
  \sup_{s\in\Omega_{\tilde\rho}}\norm{f(s)}_{L^2_\gamma}\cdot
  \set{\int_{T_0+\tilde\rho}^{T_1-\tilde\rho}\com{G_\sigma(f(s))}_{\alpha}^2ds}^{1\over2}\set{\int_{T_0+\tilde\rho}^{T_1-\tilde\rho}
  \norm{\nabla_v\p^{\alpha-1}f(s)}_{L^2_\gamma}^2ds}^{1\over2}\\
  &+&C_2\sum_{2\leq\abs\beta\leq\abs\alpha}
  C_{\alpha}^\beta
  \com{G_\sigma(f)}_{\tilde\rho,~\beta-2}\set{\int_{T_0+\tilde\rho}^{T_1-\tilde\rho}\norm{\nabla_v\p^{\alpha-\beta+1}f(s)}_{L^2_\gamma}^2ds}^{{1\over2}}\\
  &\indent\times&\set{\int_{T_0+\tilde\rho}^{T_1-\tilde\rho}\norm{\nabla_v\p^{\alpha-1}f(s)}_{L^2_\gamma}^2ds}^{{1\over2}}\\
  &+&C_2\sum_{1\leq\abs\beta\leq\abs\alpha}
  C_{\alpha}^\beta
  \com{G_\sigma(f)}_{\tilde\rho, \alpha-\beta}\set{\int_{T_0+\tilde\rho}^{T_1-\tilde\rho}\norm{\p^{\beta}f(s)}_{L^2_\gamma}^2ds}^{1\over2}
  \set{\int_{T_0+\tilde\rho}^{T_1-\tilde\rho}\norm{\nabla_v\p^{\alpha-1}f(s)}_{L^2_\gamma}^2ds}^{1\over2}\\
  &=&(S_1)+(S_2)+(S_3)+(S_4)+(S_5).
\end{eqnarray*}
In order to treat the above five terms, we need the following
estimates which are deduced directly from the the induction
hypothesis. By the truth of $(Q)_{k-1},$ we have
\begin{equation}\label{ind1}
\set{\int_{T_0+\tilde\rho}^{T_1-\tilde\rho}\norm{\nabla_v\p^{\alpha-1}f(s)}_{L^2_\gamma}^2ds}^{1/2}
 \leq \frac{A^{\abs\alpha-1}}{ \tilde\rho^{\abs\alpha-1}}\com{(\abs\alpha-2)!}^\sigma;
\end{equation}
\begin{equation}\label{ind2}
\begin{array}{lll}
&\set{\int_{T_0+\tilde\rho}^{T_1-\tilde\rho}\norm{\nabla_v\p^{\alpha-\beta+1}f(s)}_{L^2_\gamma}^2ds}^{1\over2}
\leq
\frac{A^{\abs\alpha-\abs\beta+1}}{\tilde\rho^{\abs\alpha-\abs\beta+1}}\com{(\abs\alpha-\abs\beta)!}^\sigma,\:
2\leq\abs\beta\leq\abs\alpha;
\end{array}
\end{equation}
and
\begin{equation}\label{ind3}
\set{\int_{T_0+\tilde\rho}^{T_1-\tilde\rho}\norm{\p^{\beta}f(s)}_{L^2_\gamma}^2ds}^{1/2}
 \leq
\frac{ A^{\abs\beta-1}}{
\tilde\rho^{\abs\beta-1}}\com{(\abs\beta-2)!}^\sigma,\quad
 2\leq\abs\beta\leq\abs\alpha.
\end{equation}
the last inequality using the fact $\norm {\p^\beta
f}_{L^2_\gamma}\leq \norm {\nabla_v\p^{\beta-1} f}_{L^2_\gamma} $
for any $\beta$ with $2\leq\abs\beta\leq\alpha.$ Observe that there
exists a constant $\tilde B>1,$ depending only on $B$ and $\sigma,$
such that
\begin{eqnarray*}
\left\{
\begin{array}{lll}
B^{m}\inner{m!}^\sigma\leq
 \tilde B^{m}\inner{(m-1)!}^\sigma,\quad
 1\leq m\leq\abs\alpha-1,\\
 B^{\abs\alpha}\inner{\abs\alpha!}^\sigma\leq
 \tilde  B^{\abs\alpha-1}\inner{(\abs\alpha-2)!}^\sigma.
 \end{array}
\right.
\end{eqnarray*}
With \reff{ind3}, one has, by taking $A$ such that $ A\geq \tilde
B,$
\begin{eqnarray}\label{mfv3}
 \set{\int_{T_0+\tilde\rho}^{T_1-\tilde\rho}\com{G_\sigma(f(s))}_{\alpha}^2ds}^{1/2}
 \leq\frac{2A^{|\alpha|-1}}{\tilde\rho^{|\alpha|-1}}\com{(|\alpha|-2)!}^\sigma.
\end{eqnarray}
Next we shall treat the term  $\com{G_\sigma(f)}_{\tilde\rho,
\lambda}$ which equals to  $
\sup_{s\in\Omega_{\tilde\rho}}\norm{\p^{\lambda}f(s)}_{L^2}+B^{\abs\lambda}\inner{\abs\lambda!}^\sigma$
by definition. It follows from the truth of $(Q)_{k-1}$ again that
\begin{eqnarray*}
\sup_{t\in\Omega_{\tilde\rho}}\norm{\p^{\lambda}f(t)}_{L^2}
 \leq \frac{A^{\abs\lambda}}{\tilde\rho^{\abs\lambda}}\inner{(\abs\lambda-1)!}^\sigma,\quad \forall \lambda,
 \: 1\leq \abs\lambda\leq k-1,
\end{eqnarray*}
from which and the fact $ A\geq \tilde B,$  we get the estimate on
$\com{G_\sigma(f)}_{\tilde\rho, \lambda},$ that is,
\begin{eqnarray}\label{mfv1}
 \com{G_\sigma(f)}_{\tilde\rho,\lambda}
 \leq \frac{2A^{\abs\lambda}}{\tilde\rho^{\abs\lambda}}\inner{(\abs\lambda-1)!}^\sigma, \quad 1\leq \abs\lambda\leq\abs\alpha-1.
\end{eqnarray}

Now the above estimates allow us to deal with the terms from $(S_1)$
to $(S_5).$ Note that $\norm{\p^\alpha f(s)}_{L^2}\leq
\norm{\p^\alpha f(s)}_{L^2_\gamma}$  and hence from \reff{vpi} and
\reff{ind3}, it follows immediately that
 \begin{eqnarray*}
 (S_1)\leq
 {C_3k\over\rho}\set{\frac{A^{\abs\alpha-1}}{\tilde\rho^{\abs\alpha-1}}\com{(\abs\alpha-2)!}^\sigma}^2.
\end{eqnarray*}
This along with the facts $k=\abs\alpha$ and
$\rho^{-1}<\tilde\rho^{-1}<\tilde\rho^{-2}$ shows at once
 \begin{eqnarray}\label{S1}
 (S_1)\leq C_4\set{\frac{A^{\abs\alpha-1}}{\tilde\rho^{\abs\alpha}}\com{(\abs\alpha-1)!}^\sigma}^2.
\end{eqnarray}
And by virtue of \reff{ind1}, we obtain
 \begin{eqnarray}\label{S2}
 (S_2)&\leq
  C_2\abs\alpha^2\set{\frac{A^{\abs\alpha-1}}{\tilde\rho^{\abs\alpha-1}}\com{(\abs\alpha-2)!}^\sigma}^2
   \leq C_5\set{\frac{A^{\abs\alpha-1}}{\tilde\rho^{\abs\alpha-1}}\com{(\abs\alpha-1)!}^\sigma}^2 .
\end{eqnarray}
For the term $(S_3),$ Combining \reff{A}, \reff{ind1} and
\reff{mfv3}, one has
\begin{eqnarray}\label{S4}
  (S_3)\leq
   C_6 A
  \set{\frac{A^{\abs\alpha-1}}{\tilde\rho^{\abs\alpha-1}}\com{(\abs\alpha-2)!}^\sigma}^2.
\end{eqnarray}
The treatment of the terms $(S_4)$ is a little more complicated.
Write $(S_4)=(S_4)'+(S_4)^{''}$ with
\begin{eqnarray*}
  (S_4)'&=C_2\sum\limits_{\abs\beta=2}
  C_{\alpha}^\beta
  \com{G_\sigma(f)}_{\tilde\rho, 0}\set{\int_{T_0+\tilde\rho}^{T_1-\tilde\rho}\norm{\nabla_v\p^{\alpha-\beta+1}f(s)}_{L^2_\gamma}^2ds}^{{1\over2}}
  \set{\int_{T_0+\tilde\rho}^{T_1-\tilde\rho}\norm{\nabla_v\p^{\alpha-1}f(s)}_{L^2_\gamma}^2ds}^{{1\over2}}
\end{eqnarray*}
and
\begin{eqnarray*}
  (S_4)^{''}&=&C_2\sum\limits_{3\leq\abs\beta\leq\abs\alpha}
  C_{\alpha}^\beta
  \com{G_\sigma(f)}_{\tilde\rho,\beta-2}\set{\int_{T_0+\tilde\rho}^{T_1-\tilde\rho}\norm{\nabla_v\p^{\alpha-\beta+1}f(s)}_{L^2_\gamma}^2ds}^{{1\over2}}\\
  &\indent\times&\set{\int_{T_0+\tilde\rho}^{T_1-\tilde\rho}\norm{\nabla_v\p^{\alpha-1}f(s)}_{L^2_\gamma}^2ds}^{{1\over2}}.
\end{eqnarray*}
It is easy to verify that, by \reff{A}, \reff{ind1} and \reff{ind2},
\begin{eqnarray*}
  (S_4)'\leq C_7 A \abs\alpha^2\set{
  \frac{A^{\abs{\alpha}-1}}{\tilde\rho^{\abs{\alpha}-1}}\com{(\abs\alpha-2)!}^\sigma}^2
  \leq  C_8 A\set{
  \frac{A^{\abs{\alpha}-1}}{\tilde\rho^{\abs{\alpha}-1}}\com{(\abs\alpha-1)!}^\sigma}^2.
\end{eqnarray*}
And by virtue of \reff{ind1} and \reff{mfv1}, we know $(S_4)^{''}$
is bounded from above by
\begin{eqnarray*}
 &\sum_{3\leq\abs\beta\leq\abs\alpha}
 \frac{C_3 \abs\alpha
 !}{\abs\beta!(\abs{\alpha}-\abs\beta)!}\frac{A^{|\beta|-2}}{\tilde\rho^{|\beta|-2}}\com{(|\beta|-3)!}^\sigma
 \frac{A^{\abs{\alpha}-|\beta|+1}}{\tilde\rho^{\abs{\alpha}-|\beta|+1}}\com{(\abs\alpha-|\beta|)!}^\sigma
 \frac{A^{|\alpha|-1}}{\tilde\rho^{|\alpha|-1}}\com{(|\alpha|-2)!}^\sigma,
\end{eqnarray*}
from which we get
\begin{eqnarray*}
   (S_4)^{''}\leq
  C_{9}\set{\frac{A^{\abs{\alpha}-1}}{\tilde\rho^{\abs{\alpha}-1}}}^2\com{(\abs\alpha-2)!}^\sigma
  \sum_{3\leq\abs\beta\leq\abs\alpha}
 \frac{\abs\alpha !}{\abs\beta!(\abs{\alpha}-\abs\beta)!}\com{(|\beta|-3)!}^\sigma
 \com{(\abs\alpha-|\beta|)!}^\sigma.
\end{eqnarray*}
Direct verification shows that
\begin{eqnarray*}
 \sum_{3\leq\abs\beta\leq\abs\alpha}
 \frac{\abs\alpha !}{\abs\beta!(\abs{\alpha}-\abs\beta)!}\com{(|\beta|-3)!}^\sigma
 \com{(\abs\alpha-|\beta|)!}^\sigma
 &\leq&\com{(\abs\alpha-1)!}^\sigma
 \sum_{3\leq\abs\beta\leq\abs\alpha}\frac{6\abs\alpha}{\abs\beta^3}\\
 &\leq& C_{10}
 \com{(\abs\alpha-1)!}^\sigma\abs\alpha^\sigma \quad({\rm by
 \reff{sum1}}).
\end{eqnarray*}
Combining these, we get the estimate of $(S_4)^{''}$; that is
\begin{eqnarray*}\label{S3}
 (S_4)^{''}\leq C_{11}  \set{
  \frac{A^{\abs{\alpha}-1}}{\tilde\rho^{\abs{\alpha}-1}}\com{(\abs\alpha-1)!}^\sigma}^2.
\end{eqnarray*}
This with the estimate of $(S_4)'$ gives
\begin{equation}\label{S3}
 (S_4)=(S_4)'+(S_4)^{''}\leq C_{12} A \set{
  \frac{A^{\abs{\alpha}-1}}{\tilde\rho^{\abs{\alpha}-1}}\com{(\abs\alpha-1)!}^\sigma}^2.
\end{equation}
The term $(S_5)$ can be handled exactly as above, and we have, by
virtue of \reff{sum2}, \reff{A}, \reff{ind1} and \reff{ind3},
\begin{equation}\label{S5}
 (S_5)\leq C_{13} A \set{
  \frac{A^{\abs{\alpha}-1}}{\tilde\rho^{\abs{\alpha}-1}}\com{(\abs\alpha-1)!}^\sigma}^2.
\end{equation}
Now combination of  \reff{S1}, \reff{S2}, \reff{S4},\reff{S3} and
\reff{S5}
 gives that, for all $t\in\Omega\rho,$
\begin{eqnarray*}
  \norm{\p^\alpha
  f(t)}_{L^2}^2\leq (S_1)+(S_2)+(S_3)+(S_4)+(S_5)
  \leq C_{14}A
  \set{\frac{A^{|\alpha|-1}}{\tilde\rho^{|\alpha|}}\com{(|\alpha|-1)!}^{\sigma}}^2 .
\end{eqnarray*}
Note that $\tilde\rho^{-\abs\alpha}=\tilde\rho^{-k}\leq
5\rho^{-\abs\alpha}$ by \reff{fact}, and hence the above inequality
yields
\begin{eqnarray*}
  \norm{\p^\alpha
  f(t)}_{L^2}^2
  \leq C_{15}A
  \set{\frac{A^{|\alpha|-1}}{\rho^{|\alpha|}}\com{(|\alpha|-1)!}^{\sigma}}^2,\quad
  \forall~t\in\Omega\rho .
\end{eqnarray*}
Taking $A$ large enough such that $A\geq
16\max\bigset{\sup\limits_{s\in[T_0,T_1]}\norm
{f(s)}_{H^2_\gamma}+1, \:\tilde B,\: C_{15}},$ then we get finally
\[\norm{\p^\alpha
  f(t)}_{L^2}^2\leq \set{{1\over2}\frac{A^{\abs\alpha}}{\rho^{\abs\alpha}}\com{(\abs\alpha-1)!}^\sigma}^{2}.
\]
The above inequality holds for all $t\in\Omega_\rho$, and hence
\reff{first term} follows.

In order to finish the proof, it remains to prove
\begin{equation}\label{second term}
 \set{\int_ {T_0+\rho}^{T_1-\rho}\|\nabla_v
 \p^{\tilde\alpha}f(t)\|_{L^2_{\gamma}}^2dt}^{1\over2}\leq
 {1\over2}\frac{A^{\abs\alpha}}{\rho^{\abs\alpha}}\com{(\abs\alpha-1)!}^\sigma,
\quad \forall~\abs{\tilde\alpha}=k,\quad
\forall~0<\rho<\frac{T_1-T_0}{2}.
\end{equation}
And it can be handled exactly as \reff{first term}.  The only
difference is that the multi-index $\alpha$ and the term
$\norm{\p^\alpha
  f(t)}_{L^2}^2
  =\vpi_{\rho,k}(t)\norm{\p^\alpha
  f(t)}_{L^2}^2
  -\vpi_{\rho,k}(T_0)\norm{\p^\alpha
  f(T_0)}_{L^2}^2 $ appearing in the above
argument will be replace respectively by $\tilde\alpha$ and
\[
  C_1\int_{T_0+\rho}^{T_1-\rho}\norm{\nabla_v\p^{\tilde\alpha}f(s)}_{L^2_{\gamma}}^2ds
  \leq C_1\int_{T_0+\tilde\rho}^{T_1-\tilde\rho}\vpi_{\rho, k}(s)\norm{\nabla_v\p^{\tilde\alpha}f(s)}_{L^2_{\gamma}}^2ds.
\]

Finally, combination of \reff{first term} and \reff{second term}
gives the truth of $(Q)_k.$ This completes the proof of Proposition
\ref{prp}.

\section{Proof of Lemma \ref{crucial}}
\label{section3}

For simplicity , in this section, $\p_{v_iv_j}\bar a_{ij}$ and $\bar
a_{ij}$ will stand for $\sum_{1\leq i,j\leq 3}\p_{v_iv_j}\bar
a_{ij}$ and $\sum_{i,j}\bar a_{ij},$ respectively. In the sequel $C$
is used to denote different constants which  can be replaced by a
larger one and depend only on  $\gamma,$ the Gevrey index $\sigma,$
and the initial mass $M_0$, the initial energy $E_0$ and the initial
entropy $H_0$.

Our starting point is the following uniformly ellipticity property
of the matrix $(\bar a_{ij})$, cf. Proposition 4 of
\cite{desv-vill1}.

\begin{lem}
 There exists a constant $K$, depending only on $\gamma$ and $M_0, E_0, H_0$, such
 that
 \begin{equation}\label{ellipticity}
  \sum_{1\leq i,j\leq 3}\bar a_{ij}(t,v)\xi_i\xi_j\geq K(1+\abs
  v^2)^{\gamma/2}\abs\xi^2, \quad \forall~\xi\in{\Real}^3.
 \end{equation}
\end{lem}

\begin{rem}
Although the ellipticity of $(a_{ij})$ was proved in
\cite{desv-vill1} under the hard potential case $\gamma>0,$   it's
still true for $\gamma=0,$ the Maxwellian case. This can be seen in
the proof of Proposition 4 of \cite{desv-vill1}.
\end{rem}

\begin{lem}\label{a}

There exists a constant $B$, depending only on the Gevrey index
$\sigma>1,$ such that for all multi-indices $\beta$ with
$|\beta|\geq 2$ and all $g,h\in L^2_\gamma({\Real}^3),$
\begin{eqnarray*}
&\int_{{\Real}^3}(\p^{\beta}\bar a_{ij}(t,v))g(v)h(v)dv \leq
 C\norm{g}_{L^2_{\gamma}}\norm{h}_{L^2_{\gamma}}\com{G_\sigma(f(t))}_{\beta-2},
 \quad \forall\:t>0,
\end{eqnarray*}
where $\com{G_\sigma(f(t))}_{\beta-2}=\set{
 \norm{\p^{\beta-2}f(t)}_{L^2}+B^{\abs{\beta}-2}\com{(|\beta|-2)!}^\sigma}.$
\end{lem}

{\bf Proof~} For $\sigma>1,$ there exists a function $\psi\in
G_0^\sigma({\Real}^3)$ compact support in $\set{v\in{\Real}^3\V \abs
v\leq 2},$ satisfying  $\psi(v)=1$ on the ball $\set{v\in{\Real}^3\V
\abs v\leq 1}$, and  moreover for some constant $L>4$ depending only
on $\sigma,$
\begin{equation}\label{psi}
   \sup \abs {\p^{\lambda} \psi}\leq
   L^{|\lambda|}(|\lambda|!)^\sigma,\quad \forall~\lambda.
\end{equation}
For the construction of $\psi$, see \cite{Ro} for example. Write
$a_{ij}=\psi a_{ij}+(1-\psi) a_{ij}.$ Then $\bar a_{ij}=(\psi
a_{ij})*f+[(1-\psi) a_{ij}]*f$, and hence
\[\p^{\beta}\bar
a_{ij}=\com{\p^{\tilde\beta}(\psi
a_{ij})}*(\p^{\beta-\tilde\beta}f)+\set{\p^{\beta}\com{(1-\psi)
a_{ij}}}*f,
\]
where $\tilde\beta$ is an arbitrary multi-index satisfying
$\tilde\beta\leq\beta$ and $|\tilde\beta|=2.$ We firstly treat the
term $\com{\p^{\tilde\beta}(\psi
a_{ij})}*(\p^{\beta-\tilde\beta}f)$. It is easy to verify that for
all $\tilde\beta$ with $\bigabs{\tilde\beta}=2,$
\[\bigabs{(\p^{\tilde\beta}a_{ij})(v-\vstar)}\leq C\abs{v-\vstar}^\gamma,\]
from which, we can compute
\begin{eqnarray*}
  \abs{\com{\p^{\tilde\beta}(\psi
  a_{ij})}*(\p^{\beta-\tilde\beta}f)(v)}&=\abs{\int_{{\Real}^3}
  \com{\p^{\tilde\beta}(\psi
  a_{ij})}(v-\vstar)\cdot(\p^{\beta-\tilde\beta}f)(\vstar)d\vstar}\\
  &\leq C\int_{\set{\abs{\vstar-v}\leq2}}
  \abs{(\p^{\beta-\tilde\beta}f)(\vstar)}d\vstar\\
  &\leq C\norm{\p^{\beta-2}f(t)}_{L^2}.
\end{eqnarray*}
Next  the term $\set{\p^{\beta}\com{(1-\psi) a_{ij}}}*f$, we use
Leibniz's formula to get
\begin{eqnarray*}
  &\abs{\biginner{\p^{\beta}\com{(1-\psi) a_{ij}}}*f(v)}\\&=
  \big|\sum_{0\leq\abs{\lambda}\leq\abs{\beta}}C_\beta^\lambda\int_{{\Real}^3}
  \com{\p^{\beta-\lambda}(1-\psi)}(v-\vstar)\cdot\inner{\p^{\lambda}
  a_{ij}}(v-\vstar)\cdot f(t,\vstar)d\vstar\big|
  \\&\leq \big|\sum_{0\leq\abs{\lambda}\leq\abs{\beta}-1}C_\beta^\lambda
  \int_{\set{1\leq\abs{\vstar-v}\leq2}}
  \inner{\p^{\beta-\lambda}\psi}(v-\vstar)\cdot\inner{\p^{\lambda}
  a_{ij}}(v-\vstar)\cdot f(t,\vstar)d\vstar\big|\\
  &\indent+
  \big|\int_{\set{\abs{\vstar-v}\geq1}}
  \com{(1-\psi)}(v-\vstar)\cdot\inner{\p^{\beta}
  a_{ij}}(v-\vstar)\cdot f(t,\vstar)d\vstar\big|
  \\
  &=J_1+J_2.
\end{eqnarray*}
In view of \reff{coe}, we can find a constant $\tilde C,$  such that
for all multi-indices $\lambda,$
\[
\abs{\inner{\p^{\lambda} a_{ij}}(v-\vstar)}\leq \tilde
C^{\abs\lambda}\abs\lambda! \quad {\rm
for}~1\leq\abs{\vstar-v}\leq2.
\]
And for all $\beta$ with $\abs\beta\geq2,$
\[
\abs{\inner{\p^{\beta} a_{ij}}(v-\vstar)}\leq \tilde
C^{\abs\beta}\abs\beta!(1+\abs\vstar^\gamma+\abs v^\gamma) \quad
{\rm for~} \abs{\vstar-v}\geq1.
\]
These along with \reff{psi} give the upper bound of   $J_1$ and
$J_2$,
\[
 J_1\leq L^{\abs\beta}(\abs\beta!)^\sigma\cdot\norm{f(t)}_{L^1}
 \sum_{0\leq\abs{\lambda}\leq\abs{\beta}-1}\inner{\frac{\tilde
 C}{L}}^{\abs\lambda};\quad J_2\leq2\tilde
C^{\abs\beta}(\abs\beta!)^\sigma\cdot\norm{f(t)}_{L_\gamma^1}(1+\abs
v^\gamma).
 \]
By taking $L$ large enough such that $L\geq 2\tilde C,$ then we get
\begin{eqnarray*}
  J_1+J_2\leq C\norm{f(t)}_{L_\gamma^1}
  L^{\abs\beta}(\abs\beta!)^\sigma(1+\abs v^2)^{\gamma/2}.
\end{eqnarray*}
This along with the  fact $\norm{f(t)}_{L_\gamma^1}\leq M_0+2E_0$
gives at once
\begin{eqnarray*}
 \abs{\biginner{\p_v^{\beta}\com{(1-\psi) a_{ij}}}*f(v)}\leq J_1+J_2\leq
 C L^{\abs\beta}(\abs\beta!)^\sigma(1+\abs v^2)^{\gamma/2}.
\end{eqnarray*}
 Now we choose a constant $B$ such that $L^{\abs\beta}(\abs\beta!)^\sigma\leq
B^{\abs\beta}\com{(\abs\beta-2)!}^\sigma.$ From the above inequality
it follows immediately that
\begin{eqnarray*}
  \abs{\biginner{\p_v^{\beta}\com{(1-\psi) a_{ij}}}*f(v)}
  \leq C B^{\abs\beta-2}\com{(\abs\beta-2)!}^\sigma(1+\abs v^2)^{\gamma/2}.
\end{eqnarray*}
Combining the estimate on the term $\com{\p^{2}(\psi
a_{ij})}*(\p^{\abs\beta-2}f)$,  one has
\begin{eqnarray*}
  \abs{\p^{\beta}\bar a_{ij}(v)}&\leq C
  \set{\norm{\p^{{\beta}-2}f(t)}_{L^2}
  +B^{\abs{\beta}-2}\com{(|\beta|-2)!}^\sigma\cdot(1+\abs v^2)^{\gamma/2}}\\
  &\leq C \com{G_\sigma(f(t))}_{{\beta}-2}\cdot
  (1+\abs v^2)^{\gamma/2}.
\end{eqnarray*}
Together with Cauchy's inequality, we get the desired inequality.

In quite similar argument,  we have the following
\begin{lem}\label{c}
For all multi-indices $\beta$ with $|\beta|\geq0$ and all $g,h\in
L^2_\gamma({\Real}^3),$  one has
\begin{eqnarray*}
&\int_{{\Real}^3}(\p^{\beta}\bar c(t,v))g(v)h(v)dv \leq
 C\norm{g}_{L^2_\gamma}\norm{h}_{L^2_\gamma}
 \cdot\com{G_\sigma(f(t))}_{\beta}, \quad
\forall\: t\geq0.
\end{eqnarray*}
\end{lem}

\bigskip

The rest of the paper is devoted to

\smallskip

{\em{Proof of Lemma \ref{crucial}.}} Let $b_j=\p_{v_i}
a_{ij}(v)=-2\abs v^\gamma v_j.$ Then we have the following relation
\[ \p_{v_i}\bar a_{ij}(v)=\bar b_j(v),\quad \p_{v_j}\bar b_j=\bar c.\]
Since $f$ satisfies $\p_tf=\bar a_{ij}\p_{v_ivj}f-\bar c f,$ then
direct verification shows
\begin{eqnarray*}
 \frac{d}{dt}\norm{\p^\mu f(t)}_{L^2}^2&=2\int_{{\Real}^3}\com{\p_t
 \p^\mu f(t,v)}\cdot\com{\p^\mu  f(t,v)}dv\\&=2\int_{{\Real}^3}\com{
 \p^\mu  \biginner{ \bar a_{ij}\p_{v_iv_j}f-\bar c f}}\cdot\com{\p^\mu  f(t,v)}
 dv.
\end{eqnarray*}
Moreover, using Leibniz's formula on the term $\p^\mu  \inner{ \bar
a_{ij}\p_{v_iv_j}f-\bar c f}$ , one has
\begin{eqnarray*}
  \p_t\norm{\p^\mu f(t)}_{L^2}^2&=&2\int_{{\Real}^3}
  \bar a_{ij}\inner{\p_{v_iv_j}\p^\mu  f}\cdot\inner{\p^\mu  f } dv\\
  &&+2\sum_{\abs\beta=1}
  C_{\mu }^\beta\int_{{\Real}^3}\inner{\p^\beta\bar a_{ij}}\inner{\p_{v_iv_j}\p^{\mu -\beta} f}
  \cdot\inner{\p^\mu f } dv\\
  &&+2\sum_{2\leq\abs\beta\leq\abs\mu}
  C_{\mu }^\beta\int_{{\Real}^3}\inner{\p^\beta\bar a_{ij}}\inner{\p_{v_iv_j}\p^{\mu -\beta} f}
  \cdot\inner{\p^\mu f } dv\\
  &&-2\sum_{0\leq\abs\beta\leq\abs\mu }
  C_{\mu }^\beta\int_{{\Real}^3}\inner{\p^{\mu-\beta}\bar c}\inner{
  \p^{\beta} f}\cdot\inner{\p^\mu  f }dv\\
  &=&(I)+(II)+(III)+(IV).
\end{eqnarray*}
We shall proceed to treat the above terms by the following steps.

\smallskip
 {\bf Step 1. Upper bound for the term $(I).$}

\smallskip
Integrating by parts, one has
\begin{eqnarray*}
  (I)&=&-2\int_{{\Real}^3}
  \bar a_{ij}\inner{\p_{v_j}\p^\mu f}\cdot\inner{\p_{v_i}\p^\mu f } dv-2\int_{{\Real}^3}
  \bar b_{j}\inner{\p_{v_j}\p^\mu f}\cdot\inner{\p^\mu f } dv\\
  &=&(I)_1+(I)_2.
\end{eqnarray*}
The ellipticity property \reff{ellipticity} of $(a_{ij})$ gives that
\begin{eqnarray*}
  (I)_1\leq -2K\int_{{\Real}^3}
  \abs{\nabla_v\p^\mu f}^2(1+\abs v^2)^{\gamma/2}dv
  = -2K\norm{\nabla_v\p^\mu f(t)}_{L^2_\gamma}^2.
\end{eqnarray*}
For the term $(I)_2,$ integrating by parts again, we have
\begin{eqnarray*}
  (I)_2=-(I)_2+2\int_{{\Real}^3}
  \bar c\inner{\p^\mu f}\cdot\inner{\p^\mu f } dv.
\end{eqnarray*}
This along with the fact $
  \abs{\bar c(v)}\leq C\norm{f(t)}_{L^1_\gamma}
  (1+\abs v^2)^{\gamma/2}\leq C
  (1+\abs v^2)^{\gamma/2}
$ shows immediately
\[
  (I)_2\leq C\norm{\p^\mu f(t)}_{L^2_\gamma}^2
  \leq C\norm{\nabla_v\p^{\mu-1} f(t)}_{L^2_\gamma}^2.
\]
Combining these, we get the estimate on the term $(I)$, that is
\begin{equation}\label{I}
  (I)\leq -2K\norm{\nabla_v\p^\mu f(t)}_{L^2_{\gamma}}^2+C\norm{\nabla_v\p^{\mu-1} f(t)}_{L^2_{\gamma}}^2.
\end{equation}

\smallskip
{\bf Step 2. Upper bound for the term $(II).$}
\smallskip

Recall $(II)=2\sum_{\abs\beta=1} C_{\mu
}^\beta\int_{{\Real}^3}\inner{\p^\beta\bar
a_{ij}}\inner{\p_{v_iv_j}\p^{\mu -\beta} f} \cdot\inner{\p^\mu f
}dv.$ Integrating by parts, we get
\begin{eqnarray*}
 (II)&=-2\sum_{\abs\beta=1}
  C_{\mu}^\beta\int_{{\Real}^3}\inner{\p^\beta\bar b_{j}}\inner{\p_{v_j}\p^{\mu-\beta} f}
  \cdot\inner{\p^\mu f } dv\\
  &\indent-2\sum_{\abs\beta=1}
  C_{\mu}^\beta\int_{{\Real}^3}\inner{\p^\beta\bar a_{ij}}\inner{\p_{v_j}\p^{\mu-\beta} f}
  \cdot\inner{\p_{v_i}\p^\mu f } dv\\
  &=(II)_1+(II)_2.
\end{eqnarray*}
Note $\abs{\p^\beta\bar b_{j}(t,v)}\leq C(1+\abs v^2)^{\gamma/2}$
for any $\beta$ with $\abs\beta=1$ and hence
\begin{eqnarray*}
 (II)_1\leq C\abs\mu\cdot\norm{\nabla_v\p^{\mu-\beta}f(t)}_{L^2_\gamma}
 \norm{\p^{\mu}f(t)}_{L^2_\gamma}
 \leq
 C\abs\mu\cdot\norm{\nabla_v\p^{\mu-1}f(t)}_{L^2_\gamma}^2.
\end{eqnarray*}
For the term $(II)_2$, noticing that $\mu=\beta+(\mu-\beta),$ it can
be rewritten as the following form
\begin{eqnarray*}
 (II)_2&=-2\sum_{\abs\beta=1}
  C_{\mu}^\beta\int_{{\Real}^3}\inner{\p^\beta\bar a_{ij}}\inner{\p_{v_j}\p^{\mu-\beta} f}
  \cdot\inner{\p^\beta\p_{v_i}\p^{\mu-\beta} f } dv.
\end{eqnarray*}
Since $\abs \beta=1,$ we can integrate by parts to get
\begin{eqnarray*}
  (II)_2=&2\sum_{\abs\beta=1}
  C_{\mu}^\beta\int_{{\Real}^3}\inner{\p^{\beta}\bar a_{ij}}\inner{\p_{v_j}\p^{\mu} f}
  \cdot\inner{\p_{v_i}\p^{\mu-\beta} f } dv\\
  &+2\sum_{\abs\beta=1}
  C_{\mu}^\beta\int_{{\Real}^3}\inner{\p^{\beta+\beta}\bar a_{ij}}\inner{\p_{v_j}\p^{\mu-\beta} f}
  \cdot\inner{\p_{v_i}\p^{\mu-\beta} f } dv\\
  =&-(II)_2+2\sum_{\abs\beta=1}
  C_{\mu}^\beta\int_{{\Real}^3}\inner{\p^{\beta+\beta}\bar a_{ij}}\inner{\p_{v_j}\p^{\mu-\beta} f}
  \cdot\inner{\p_{v_i}\p^{\mu-\beta} f } dv.
\end{eqnarray*}
Hence
\begin{eqnarray*}
  (II)_2&=\sum_{\abs\beta=1}
  C_{\mu}^\beta\int_{{\Real}^3}\inner{\p^{\beta+\beta}\bar a_{ij}}\inner{\p_{v_j}\p^{\mu-\beta} f}
  \cdot\inner{\p_{v_i}\p^{\mu-\beta} f } dv.
\end{eqnarray*}
This along with the fact $\abs{\p^{\beta+\beta}\bar a_{ij}(v)}\leq C
(1+\abs v^2)^{\gamma/2}$ for all $\beta$ with $\abs\beta=1$ shows at
once
\begin{eqnarray*}
 (II)_2\leq C\sum_{\abs\beta=1}
  C_{\mu}^\beta\cdot\norm{\nabla_v\p^{\mu-\beta}f}_{L^2_{\gamma}}^2
  \leq C\abs\mu\cdot\norm{\nabla_v\p^{\mu-1}f}_{L^2_{\gamma}}^2.
\end{eqnarray*}
Combining these, we obtain
\begin{eqnarray}\label{II}
 (II)\leq C\abs\mu\cdot\norm{\nabla_v\p^{\abs\mu-1}f}_{L^2_{\gamma}}^2.
\end{eqnarray}

\smallskip
{\bf Step 3. Upper bound for the term $(III)$ and $(IV)$ and the
conclusion.}
\smallskip

Recall $(III)=2\sum_{2\leq\abs\beta\leq\abs\mu}
  C_{\mu }^\beta\int_{{\Real}^3}\inner{\p^\beta\bar a_{ij}}\inner{\p_{v_iv_j}\p^{\mu -\beta} f}
  \cdot\inner{\p^\mu f }\,dv$ and
\[
  (IV)=-2\sum_{0\leq\abs\beta\leq\abs\mu }
  C_{\mu }^\beta\int_{{\Real}^3}\inner{\p^{\mu -\beta}\bar c}\inner{
  \p^\beta f}\cdot\inner{\p^\mu  f }dv.
\]
By virtue of Lemma \ref{a} and lemma \ref{c},  it follows that
\begin{eqnarray}\label{III}
\begin{array}{lll}
  &(III)\leq C \sum_{2\leq\abs\beta\leq\abs\mu}
  C_{\mu }^\beta \norm{\p_{v_iv_j}\p^{\mu -\beta}f(t)}_{L^2_\gamma}
  \cdot\norm{\p^\mu f(t)}_{L^2_\gamma}\com{G_\sigma(f(t))}_{\beta-2}\\
  &\leq C \sum_{2\leq\abs\beta\leq\abs\mu}
  C_{\mu }^\beta \norm{\nabla_v\p^{\mu -\beta+1} f(t)}_{L^2_\gamma}
  \cdot\norm{\nabla_v\p^{\mu-1}
  f(t)}_{L^2_\gamma}\com{G_\sigma(f(t))}_{\beta-2},
\end{array}
\end{eqnarray}
and
\begin{equation}\label{IV}
  (IV)\leq C\sum_{0\leq\abs\beta\leq\abs\mu}
  C_{\mu }^\beta\norm{
  \p^\beta f(t)}_{L^2_\gamma}\cdot\norm{\nabla_v\p^{\mu-1}  f(t) }_{L^2_\gamma}\cdot\com{G_\sigma(f(t))}_{{\mu}-\beta}.
\end{equation}
Combination of \reff{I}-\reff{IV}  gives the desired inequality
\begin{eqnarray*}
  \frac{d}{dt}&\norm{\p^\mu f(t)}_{L^2}^2+C_1\norm{\nabla_v\p^\mu  f(t) }_{L^2_\gamma}^2
  \leq C_2\abs\mu^2\norm{\nabla_v\p^{\mu-1} f(t)}_{L^2_\gamma}^2\\
  &+C_2 \sum_{2\leq\abs\beta\leq\abs\mu}
  C_{\mu }^\beta \norm{\nabla_v\p^{\mu -\beta+1} f(t)}_{L^2_\gamma}
  \cdot\norm{\nabla_v\p^{\mu-1}
  f(t)}_{L^2_\gamma}\cdot\com{G_\sigma(f(t))}_{\beta-2}\\
  &+C_2\sum_{0\leq\abs\beta\leq\abs\mu}
   C_{\mu }^\beta\norm{
  \p^\beta f(t)}_{L^2_\gamma}\cdot\norm{\nabla_v\p^{\mu-1}  f(t)
  }_{L^2_\gamma}\cdot\com{G_\sigma(f(t))}_{{\mu}-\beta},
\end{eqnarray*}
where $C_1, C_2$ are two constants depending only on $M_0, E_0, H_0,
\sigma$ and $\gamma.$  This completes the proof of Lemma
\ref{crucial}.

\end{document}